%% file: article/main.tex
\documentclass[9pt]{amsart}

\input{article/preamble/packages}
\input{article/preamble/preamble_formatting}
\input{article/preamble/commands}

\begin{document}

\input{article/frontmatter/title_page}

\input{article/content/abstract}

\maketitle

\input{article/content2/introduction_2}

\section*{Acknowledgements}

The author would especially like to thank Harry Schmidt for helpful conversations and guidance, without whom this paper would not exist. Additional thanks are owed to Sam Chow and Akshat Mudgal for introducing us to the problem, to Nathan Lockwood, and to Gareth Jones. We also thank Philipp Habegger for motivating us to extend Theorem \ref{additive_relations_in_irrational_powers} from the $r = s = 2$ case, established in a previous version of this article, to the general case. Joseph Harrison is supported by the Warwick Mathematics Institute Centre for Doctoral Training, and gratefully acknowledges the funding from the UK Engineering and Physical Sciences Research Council (Grant number: \texttt{EP/W524645/1}).



\input{article/content2/preliminaries}
\input{article/content2/functional_transcendence}
\input{article/content2/proof_of_theorem_B}
\input{article/content2/proof_of_theorem_C}
\input{article/content2/rational_exponents}

\clearpage

\input{article/backmatter/bibliography}

\end{document}

%% file: article/preamble/packages.tex
\usepackage{graphicx}
\usepackage{amsmath}
\usepackage{amssymb}
\usepackage{amsfonts}
\usepackage{mathrsfs}
\usepackage{enumerate}
\usepackage{lipsum}
\usepackage{nomencl}
\usepackage{tikz-cd}
\usepackage[UKenglish]{babel}
\usepackage{xcolor}
\usepackage[a4paper]{geometry}
\usepackage{bm}

%% file: article/preamble/preamble_formatting.tex
\usepackage{hyperref}
\hypersetup{
colorlinks=true, 
linkcolor=red, 
urlcolor=magenta, 
citecolor=blue 
}

\usepackage{amsthm}
\newtheoremstyle{my_theorem_style}
{}
{}
{\it}
{}
{\bfseries}
{.}
{ }
{}
\theoremstyle{my_theorem_style}
\newtheorem{definition}{Definition}
\numberwithin{definition}{section}

\newtheorem{theorem}[definition]{Theorem}
\newtheorem{lemma}[definition]{Lemma}

\newtheorem{corollary}[definition]{Corollary}

\newtheorem*{corollary*}{Corollary}
\renewenvironment{proof}{\noindent {\textbf{Proof.}}}{\\\textbf{Q.E.D.}}

\newtheorem{manualtheoreminner}{Theorem}
\newenvironment{manualtheorem}[1]{%
  \IfBlankTF{#1}
    {\renewcommand{\themanualtheoreminner}{\unskip}}
    {\renewcommand\themanualtheoreminner{#1}}%
  \manualtheoreminner
}{\endmanualtheoreminner}



%% file: article/preamble/commands.tex
\newcommand{\abs}[1]{\left|#1\right|}

\newcommand{\C}{\mathbb{C}}

\newcommand{\eps}{\varepsilon}

\newcommand{\G}{\mathbb{G}}

\renewcommand{\P}{\mathbb{P}}
\renewcommand{\phi}{\varphi}

\newcommand{\Q}{\mathbb{Q}}

\DeclareMathOperator{\rank}{rank}
\newcommand{\R}{\mathbb{R}}
\renewcommand{\rho}{\varrho}

\DeclareMathOperator{\trdeg}{trdeg}

\renewcommand{\vec}[1]{\bm{#1}}

\newcommand{\Z}{\mathbb{Z}}

%% file: article/frontmatter/title_page.tex

\title{Additive relations in irrational powers}
\author{Joseph Harrison}
\address{Mathematics Institute, University of Warwick, Coventry, UK}
\email{joseph.s.harrison@warwick.ac.uk}

%% file: article/content/abstract.tex



\begin{abstract}
    We investigate the interaction between raising to an irrational power and addition of real numbers. Thus, for a finite set $A$ of non-negative real numbers, let $A^{[c]} = \{a^c : a \in A\}$. When $k$ is a positive integer, $c$ is a real irrational number, and $A$ is a subset of an $N$-term arithmetic progression in $\R_{\geq 0}$ having cardinality at least a power of $\log{N}$, we prove that the $k$-fold sumset $|kA^{[c]}| \sim_k |A|^k/k!$ as $|A| \to \infty$. This result is uniform in $c$. When $A = \{1, \dots, N\}$ and $k = 2$, this result can be combined with existing works to show that $|A^{[c]} + A^{[c]}| \sim N^2/2$ as $N \to \infty$ whenever $c \in \R \setminus \{0, 1, 2\}$. The sumset lower bound follows from a bound on the number of equal sums of $r$ and $s \geq r$ elements of $A^{[c]}$ (by taking $r = s = k$). When $r = s = 2$ or $s > r$, our bound is optimal up to a power of $\log N$. This bound is proved using a functional transcendence theorem for certain endomorphisms of $\R_{>0}^n$, and innovations in the Pila--Wilkie counting theorem in $\R_{\exp}$ due to Binyamini, Novikov and Zak.

    In a different direction, we provide a Diophantine approximation criterion on $c$ that, when satisfied, ensures that a linear form in the $c$-th powers of multiplicatively independent integers does not vanish. The proof involves linear forms in logarithms. This provides a new proof of a fact, due to Bays--Kirby--Wilkie and Jones--Servi, that when $A$ is a multiplicatively independent set of positive integers, there are infinitely many effectively computable real numbers $c$ such that $A^{[c]}$ is linearly independent over $\Q$.
\end{abstract}

%% file: article/content2/introduction_2.tex
\section{Introduction}

Let $A \subseteq \R_{\geq 0}$ be a finite set and let $c$ be a real number. We let $A^{[c]} = \{a^c : a \in A\}$ denote the image of $A$ under the power function sending $x$ to $x^c$. In this paper, we are interested in the additive structure of the set $A^{[c]}$. Thus we are interested in the $k$-fold sumsets
\begin{align*}
    k A^{[c]} = A^{[c]} + \dots + A^{[c]} =
    \{a_1^c + \dots + a_k^c : a_1, \dots, a_k \in A\}
\end{align*}
and the $k$-fold additive energies
\begin{align*}
    E_k(A^{[c]}) = 
    \left|
    \left\{
    (\vec{a}, \vec{b}) \in A^k \times A^k :
    a_1^c + \dots + a_k^c = 
    b_1^c + \dots + b_k^c
    \right\}
    \right|.
\end{align*}
Our main result concerns lower bounds on the cardinality of $k A^{[c]}$ when $c$ is irrational, and $A$ is assumed to lie in an arithmetic progression. The result is uniform in $c$.

\begin{theorem}
[Expansion]
\label{expansion}
    Let $k$ be a positive integer. There exists an effectively computable constant $C_1 > 1$, depending on $k$, with the following property. Let $A$ be a subset of an $N$-term arithmetic progression in $\R_{\geq 0}$, such that $|A| \geq (\log N)^{C_1}$. Let $c$ be a real irrational number. Then
    \begin{align*}
        |k A^{[c]}| \geq \frac{1}{k!}(1 - o_k(1))|A|^k
    \end{align*}
    as $|A| \to \infty$. In particular, $|k A^{[c]}| \sim_k |A|^k/k!$ as $|A| \to \infty$.
\end{theorem}

In Section \ref{notations_and_conventions}, we recall the meaning of certain asymptotic notations appearing throughout this paper.

When $A = \{1, \dots, N\}$ and $c$ is instead a positive integer, these questions have been addressed, essentially, by Hooley \cite{Hooley_on_the_representations_of_a_number_as_a_sum_of_two_cubes, Hooley_on_the_numbers_that_are_representable_as_the_sum_of_two_cubes, Hooley_on_the_representation_of_a_number_as_the_sum_of_two_h-th_powers, Hooley_on_another_sieve_method_and_the_numbers_that_are_a_sum_of_two_h-th_powers}, Skinner--Wooley \cite{Skinner--Wooley_sums_of_two_kth_powers}, Browning \cite{Browning_equal_sums_of_two_kth_powers}, Heath-Brown \cite{Heath-Brown_the_density_of_rational_points_on_curves_and_surfaces, Heath-Brown_Cayleys_cubic} and Salberger \cite{Salberger_counting_rational_points_on_projective_varieties}, and many other authors. In Section \ref{rational_exponents}, we explain how these works can be combined with a linear independence result of Carr and O'Sullivan \cite{Carr--O'Sullivan_on_the_linear_independence_of_roots} to ascertain the asymptotic behaviour of $|A^{[c]} + A^{[c]}|$ when $c$ is rational. Together with our Theorem \ref{expansion}, this proves the following corollary.

\begin{corollary}
[Expansion for all exponents]
\label{expansion_for_all_exponents}
    Let $A = \{1, \dots, N\}$ and let $c \in \R \setminus \{0, 1, 2\}$. Then
    \begin{align*}
        |A^{[c]} + A^{[c]}| \sim \frac{1}{2} N^2
    \end{align*}
    as $N \to \infty$.
\end{corollary}

It is clear that the corollary does not hold if $c \in \{0, 1\}$. A theorem of Landau \cite{Landau_on_the_partition_of_positive_integers_in_four_classes_according_to_the_minimal_number_of_squares_needed_to_their_additive_composition} states that $|A^{[2]} + A^{[2]}| \ll N^2(\log{N})^{-1/2}$, so the corollary cannot hold in this case either.

Theorem \ref{expansion} can be deduced, via the Cauchy--Schwartz inequality, from an asymptotic formula for the $k$-fold additive energy of $A^{[c]}$. This, in turn, follows from the following general result, proved in Section \ref{proof_of_theorem_B}.

\begin{theorem}
[Additive relations in irrational powers]
\label{additive_relations_in_irrational_powers}
    Let $r \leq s$ be positive integers. There exists an effectively computable constant $C_2 > 0$, depending on $s$, with the following property. Let $A$ be a subset of an $N$-term arithmetic progression in the non-negative reals. Let $c$ be a real irrational number. The number of solutions to the equation
    \begin{align}
    \label{additive_relations_in_irrational_powers_equation}
        \sum_{i = 1}^s a_i^c = \sum_{j = 1}^r b_j^c
    \end{align}
    with $a_1, \dots, a_s, b_1, \dots, b_r \in A$ and $(a_1, \dots, a_s)$ not a permutation of $(b_1, \dots, b_r)$ is
    \begin{align*}
        O_s\left(
        |A|^{\theta}(\log N)^{C_2}
        \right),
    \end{align*}
    where $\theta = \max(1, \min(r, s - 1))$.
\end{theorem}

When $r = s$, the solutions with $(a_1, \dots, a_s)$ a permutation of $(b_1, \dots, b_r)$ provide a main term $r! |A|^r$. As a special case, we have an asymptotic formula for the $k$-fold additive energy
\begin{align}
\label{k-fold_energy_asymptotic}
    E_k(A^{[c]}) = k!|A|^k + O_k(|A|^{k - 1}(\log{N})^{C_2})
\end{align}
whenever $c$ is irrational, which proves Theorem \ref{expansion} with $C_1$ any real number larger than $C_2$. 

Suppose that $k = 2$, $A = \{1, \dots, N\}$, and that the equation (\ref{additive_relations_in_irrational_powers_equation}) admits a non-trivial solution; that is, a solution with $(a_1, a_2)$ not a permutation of $(b_1, b_2)$ (we will produce examples of such $c$ below). Then (\ref{additive_relations_in_irrational_powers_equation}) admits $\gg N$ solutions given by dilates of $(a_1, a_2, b_1, b_2)$. This shows that the error term in (\ref{k-fold_energy_asymptotic}) is optimal up to the power of $\log N$ for $k = 2$ and $A = \{1, \dots, N\}$. It seems that the best result on $E_2(A^{[c]})$ up to this point was the bound $O_\eps(N^{2 + \eps})$, which follows from \cite[Theorem 2]{Robert--Sargos_three-dimensional_exponential_sums_with_monomials}.

When $r < s$ and $A = \{1, \dots, N\}$, the following construction shows that Theorem \ref{additive_relations_in_irrational_powers} is also optimal up to a power of $\log N$. Let $a_1 + \dots + a_{s - r + 1} < b_1$ be multiplicatively independent positive integers. By the intermediate value theorem, we may choose $c \in (0, 1)$ so that
\begin{align*}
    a_1^c + \dots + a_{s - r + 1}^c = b_1^c.
\end{align*}
Such a $c$ must be irrational by Galois theory. Then the equation
\begin{align*}
    x_1^c + \dots + x_{s - r + 1}^c = y^c
\end{align*}
has $\gg N$ solutions in $A$ given by dilates of the solution $(a_1, \dots, a_{s - r + 1}, b_1)$. Setting the variables $b_i = a_{s - r + i}$ for $i \in \{2, \dots, r\}$ then yields $\gg N^r$ solutions to the equation (\ref{additive_relations_in_irrational_powers_equation}).

We now consider conditions on the number $c$ that might forbid solutions to equations such as (\ref{additive_relations_in_irrational_powers_equation}). For example, it can be shown that the equation
\begin{align}
\label{energy_equation}
    n_1^c + n_2^c = n_3^c + n_4^c
\end{align}
can admit non-trivial solutions, those with $(n_1, n_2)$ not a permutation of $(n_3, n_4)$, only if
\begin{align}
\label{admissable_range}
    \frac{1}{6N^2(\log{N})} \leq
    c \leq (\log{2})N.
\end{align}
To prove the upper bound, order the variables $n_1 > n_3 \geq n_4 > n_2$ so that
\begin{align*}
    n_1^c \leq n_1^c + n_2^c = n_3^c + n_4^c \leq 2n_3^c,
\end{align*}
and use $n_1 \geq n_3 + 1$. To prove the lower bound, expand $n_i^c = \exp(c \log{n_i})$ using the power series for the exponential function to obtain
\begin{align*}
    \frac{|n_1n_2 - n_3n_4|}{N^2} \leq \log\left(\frac{n_1n_2}{n_3n_4}\right) \leq 4 c N^c (\log{N}).
\end{align*}
Thus, if $c$ is small enough, then $n_1n_2 = n_3n_4$, which is incompatible with (\ref{energy_equation}) unless the solution is trivial.

When $c$ is in this admissable range (\ref{admissable_range}), the equation (\ref{energy_equation}) can admit non-trivial solutions. Let $c_n = \log_2(\phi_n)$, where $\phi_n$ is the unique solution of the equation
\begin{align*}
    \phi_n^{2n} + 1 = \phi_n^{n + 1} + \phi_n^n
\end{align*}
with $\phi_n > 1$. Then $(2^{2n}, 1, 2^{n + 1}, 2^n)$ is a non-trivial solution to (\ref{energy_equation}) for $c = c_n$. The equation satisfied by $\phi_n$ shows that $c_n \to 0$ as $n \to \infty$, and therefore the set of real numbers $c$ such that (\ref{energy_equation}) admits non-trivial solutions is not discrete. By the Gelfond--Schneider theorem \cite{Gelfond_on_the_seventh_problem_of_D._Hilbert, Schneider_transcendence_investigations}, each $c_n$ is also transcendental.

The solutions we have constructed lie in a very multiplicatively structured set, namely a geometric progression. Contrary to the last example, we might actually expect $A^{[c]}$ to be more additively unstructured when $A$ is multiplicatively structured. This is in view of the sum-product phenomenon originating in work of Erd\H{o}s and Szemer{\' e}di \cite{Erdos--Szemeredi_on_sums_and_products_of_integers}, which asserts the existence of some $\delta > 0$ such that if $A \subseteq \C$ then either $|A \cdot A|$ or $|A + A|$ should have size at least $|A|^{1 + \delta}$. Improvements on the value of $\delta$ have been obtained by many authors, for example, Cushman \cite{Cushman_a_note_on_the_sum-product_problem_and_the_convex_sumset_problem} and Rudnev and Stevens \cite{Rudnev--Stevens_an_update_on_the_sum-product_problem}. During the preparation of this manuscript, Bloom--Sawin--Schildkraut--Zhelezov \cite{Bloom--Sawin--Schildkraut--Zhelezov} have announced a counterexample to the conjecture of Erd\H{o}s and Szemer{\' e}di that one can take $\delta$ arbitrarily close to $1$, for all finite sets $A \subseteq \R$. In particular, one cannot expect to obtain results of a similar quality as Theorem \ref{expansion} from the sum-product phenomenon alone.

The next examples, and our Theorem \ref{lower_bound_on_linear_form_in_irrational_powers}, address the other extreme, where $A$ is a multiplicatively independent subset of $\{1, \dots, N\}$. Suppose that $c$ and $\log{n_1}, \dots, \log{n_s}$ are linearly independent over $\Q$. Baker's theorem \cite{Baker_linear_forms_in_the_logarithms_of_algebraic_numbers} implies that this holds when $c$ is algebraic and $n_1, \dots, n_s$ are multiplicatively independent. A deep and influential conjecture of Schanuel would, if true, imply that $\trdeg_\Q(c, n_1^c, \dots, n_s^c) \geq s$. Thus Schanuel's conjecture implies that if $c$ is algebraic and $A$ is a multiplicatively independent set of integers, then equations such as (\ref{additive_relations_in_irrational_powers_equation}), or even polynomial equations of higher degree, cannot admit multiplicatively independent solutions.

A version of Schanuel's conjecture for certain real powers has been established unconditionally by Bays, Kirby and Wilkie \cite[Theorem 1.1]{Bays--Kirby--Wilkie_A_Schanuel_property_for_exponentially_transcendental_powers}. To state their result, we introduce the notion of an exponentially transcendental number. A real number is called exponentially algebraic if it is a coordinate in a smooth solution $(y_1, \dots, y_n)$ of a system of polynomial equations in $y_1, \dots, y_n$ and $\exp(y_1), \dots, \exp(y_n)$. An exponentially transcendental number is one that is not exponentially algebraic. In \cite{Bays--Kirby--Wilkie_A_Schanuel_property_for_exponentially_transcendental_powers}, the authors prove that if $c$ is exponentially transcendental then
\begin{align*}
    \trdeg_{\Q(c)}(y_1, \dots, y_n, y_1^c, \dots, y_n^c) \geq n
\end{align*}
for any multiplicatively independent $y_1, \dots, y_n \in \R_{>0}$. Thus if $c$ is exponentially transcendental, and $A$ is a multiplicatively independent set of integers, then equations of the form (\ref{additive_relations_in_irrational_powers_equation}) will again not admit solutions without having some variables equal.

In this theme of forbidding solutions in multiplicatively independent sets, we have the following result, proved in Section \ref{proof_of_theorem_C}.

\begin{theorem}
[Non-vanishing for well-approximated exponents]
\label{lower_bound_on_linear_form_in_irrational_powers}
    Let $-B \leq a_1, \dots, a_s \leq B$ be integers, not all zero, and let $N$ be a positive integer. There exists an effectively computable, positive-valued function $\psi(B, N, s, q)$ such that the following holds. Suppose $c$ is an irrational real number, and there exist integers $a$ and $q > 0$ with
    \begin{align}
    \label{Diophantine_approximation_condition}
        \frac{a}{q} < c \leq \frac{a}{q} + \psi(B, N, s, q).
    \end{align}
    Then the form $F(n_1, \dots, n_s) = a_1 n_1^c + \dots + a_s n_s^c$ is non-zero for all multiplicatively independent $n_1, \dots, n_s \in \{1, \dots, N\}$.
\end{theorem}

Even if the form $F$ and the integer $N$ varies, we can produce infinitely many effectively computable numbers $c$ satisfying the hypothesis of Theorem \ref{lower_bound_on_linear_form_in_irrational_powers} for all $B, N$ and $s$. This is explained in Lemma \ref{existence_of_uncountably_many_c_and_effectively_computable_examples}. In particular, if $S$ is a set of multiplicatively independent positive integers then there are infinitely many effectively computable real numbers $c$ such that the set $S^{[c]}$ is linearly independent over $\Q$. This conclusion can also be deduced by applying the transcendence result \cite[Theorem 1.1]{Bays--Kirby--Wilkie_A_Schanuel_property_for_exponentially_transcendental_powers} to the effectively computable examples of exponentially transcendental numbers constructed by Jones and Servi \cite[Section 4]{Jones--Servi_on_the_decidability_of_the_real_field_with_a_generic_power_function}.

\subsection{Outline of the paper}
\label{outline_of_the_paper}

In Section \ref{preliminaries} we recall facts about algebraic tori, o-minimal geometry, and point counting that will be used to prove Theorem \ref{additive_relations_in_irrational_powers}.

In Section \ref{functional_transcendence_section}, we will prove the following functional transcendence theorem. The notion of restricting to a morphism of algebraic groups is given in Definition \ref{restricting_to_a_morphism_of_algebraic_groups}, but let us just say that if $c_1, \dots, c_n$ are real irrational numbers, then $\phi(\vec{x}) = (x_1^{c_1}, \dots, x_n^{c_n})$ does not restrict to a morphism of algebraic groups (Lemma \ref{lemma/raising_to_irrational_powers_not_a_morhpism_of_algebraic_subgroups}).

\begin{manualtheorem}{3.4}
    [Functional transcendence]
    Let $X$ be an irreducible, semi-algebraic subset of $\R_{>0}^n$. Let $\phi : \R_{>0}^n \to \R_{>0}^n$ be a continuous group homomorphism that does not restrict to a morphism of algebraic groups. Then the Zariski closure of $\phi(X)$ in $\G_m^n$ is a translate of a connected algebraic subgroup.
\end{manualtheorem}

Theorem \ref{functional_transcendence_theorem} is proved using Ax's theorem \cite[Theorem 3]{Ax_on_Schanuel's_conjectures} for the exponential function.

In Section \ref{proof_of_theorem_B} we prove Theorem \ref{additive_relations_in_irrational_powers}. We briefly explain the argument, with special emphasis on the case $r = s = 2$, which corresponds to the asymptotic formula for $E_2(A^{[c]})$. In this case, after a dilation of $A$, it suffices to count $n_1, n_2, n_3, n_4 \in \{0, \dots, N - 1\}$ satisfying
\begin{align}
\label{r=s=2_o-minimal_counting_problem}
    (\alpha + n_1)^c + (\alpha + n_2)^c = (\alpha + n_3)^c + (\alpha + n_4)^c
\end{align}
for some $\alpha \in \R_{\geq 0}$, which is a problem of counting rational points on a set definable in the o-minimal structure $\R_{\exp}$. Such problems are addressed by the celebrated counting theorem of Pila and Wilkie \cite[Theorem 1.8]{Pila--Wilkie_the_rational_points_on_a_definable_set}, which states that if $X \subseteq \R^n$ is definable in an o-minimal structure, then the number of rational points of height at most $N$ lying on $X$, and not lying on any semi-algebraic curve in $X$, is $O_{\eps, X}(N^\eps)$ for any $\eps > 0$. The question of replacing $N^\eps$ with a power of $\log N$ for sets definable in the o-minimal structure $\R_{\exp}$ is known as Wilkie's conjecture, and was recently answered positively in a breakthrough work of Binyamini--Novikov--Zak \cite{Binyamini--Novikov--Zak_Wilkie's_conjecture_for_Pfaffian_structures}. We shall make use of their result. O-minimal point counting in the fibres of definable families is what allows us to achieve uniformity in $c$ and the $N$-term arithmetic progression containing $A$.

Unfortunately, the set defined by (\ref{r=s=2_o-minimal_counting_problem}) is covered by semi-algebraic curves. We thus employ a strategy of fixing certain variables. The resulting definable sets are no longer covered by semi-algebraic curves, but we must accept a power of $|A|$ into the final bound when the fixed variables are allowed to vary again. Theorem \ref{functional_transcendence_theorem} is used to show that the solutions lying on semi-algebraic curves solve certain two-equation systems of the form (\ref{additive_relations_in_irrational_powers_equation}), and it is here that the additive nature of (\ref{additive_relations_in_irrational_powers_equation}) combines with the multiplicative subgroups produced by Theorem \ref{functional_transcendence_theorem} to severely restrict the possible semi-algebraic curves. In the case $r = s = 2$, there are two systems of two equations expressing the condition that $(a_1, a_2)$ is a permutation of $(b_1, b_2)$, and so the semi-algebraic curves only contribute to the main term in the asymptotic (\ref{k-fold_energy_asymptotic}) for $k = 2$. In the general case, each equation involves fewer than $s$ variables, and so the contribution from the semi-algebraic curves can be handled by induction. The number $s - 1$ of fixed variables is chosen to balance the contribution from the semi-algebraic curves and the o-minimal counting.

We expect that the method should generalise to allow the variables $a_1, \dots, a_s, b_1, \dots, b_r$ to lie in different sets, each lying in possibly different arithmetic progressions, and to replace the terms $a_i^c, b_j^c$ with $a_i^{c_i}, b_j^{d_j}$ for possibly distinct real irrational numbers $c_1, \dots, c_s$ and $d_1, \dots, d_r$. The method should also allow one to treat the asymmetric additive energies such as
\begin{align*}
    n_1^c + n_2^c = n_3^c + n_4^c + k
\end{align*}
for $n_1, n_2, n_3, n_4 \in A$ and some $k \in \R$, as considered in \cite{McGrath_on_the_asymmetric_additive_energy_of_polynomials}. We thank Sam Chow for this remark.

The proof of Theorem \ref{lower_bound_on_linear_form_in_irrational_powers} is carried out in Section \ref{proof_of_theorem_C}. The idea of the proof is to write
\begin{align*}
    F(\vec{n}) = \alpha + (c - a/q) \Lambda
    + O_{B, N, s}((c - a/q)^2),
\end{align*}
where $\alpha$ is an algebraic integer and $\Lambda$ is a linear form in the logarithms of $n_1, \dots, n_s$. If $\alpha$ does not vanish, a simple lower bound on $|\alpha|$ can be used to show $F(\vec{n})$ does not vanish, and if $\alpha$ does vanish, then an effective linear forms in logarithms result of Fel'dman \cite[Theorem 1]{Fel'dman_an_improvement_of_the_estimate_of_a_linear_form_in_the_logarithms_of_algebraic_numbers} can be used. In each case, $c - a/q$ must be suitably small, and this is where the function $\psi(B, N, s, q)$ enters.

Finally, in Section \ref{rational_exponents}, we discuss previous work regarding the equation (\ref{energy_equation}) when $c$ is an integer, and explain how to deduce Corollary \ref{expansion_for_all_exponents}.

\subsection{Further questions}

We now discuss some further questions that are suggested by the main results here.

\begin{enumerate}[(1)]
    \item (GAPs of higher rank) 
    We have not been able to adapt our argument to handle subsets of arbitrary generalised arithmetic progressions (GAPs). It may be the case that one can handle GAPs in the algebraic numbers by using a version of Theorem \ref{bnz_counting_theorem} for counting algebraic points. It would be of particular interest to prove bounds with an explicit dependence on the rank $r$, which would involve proving o-minimal point counting theorems with an explicit dependence on the dimension of the definable set. If $A$ is a set of non-negative real numbers satisfying $|A + A| \leq K |A|$, then some variant of the Freiman--Rusza theorem in additive combinatorics (e.g., \cite{Green--Rusza_Freimans_theorem_in_an_arbitrary_abelian_group}) would imply that $A$ lies in a GAP of rank and cardinality bounded explicitly in terms of $K$, and $K$ and $|A|$, respectively. Thus one would obtain a bound on the number of solutions to equations such as (\ref{additive_relations_in_irrational_powers_equation}) for arbitrary finite sets $A$.
    \item (Sparser subsets) One can consider Theorem \ref{expansion} for sparser sets, e.g., sets of cardinality $\log\log N$, say. If such a result is to follow from a suitable improvement of Theorem \ref{additive_relations_in_irrational_powers}, then we would require significant improvements in the Pila--Wilkie theorem. It is possible that some argument bypassing additive energy could be employed.
    \item (Higher dimensions)
    One can ask for generalisations of Theorem \ref{expansion}, or indeed Theorem \ref{additive_relations_in_irrational_powers}, to subsets of $\R_{\geq 0}^n$, where now instead of raising to an irrational power one can take the image of $A$ under a general continuous homomorphism $\phi : \R_{>0}^n \to \R_{>0}^n$ that does not restrict to a morphism of algebraic groups, in the sense of Definition \ref{restricting_to_a_morphism_of_algebraic_groups}. This would provide the analogue of Theorem \ref{additive_relations_in_irrational_powers} for systems of equations.
    \item (Diophantine inequalities)
    It would be beneficial to extend our o-minimal point counting methods to handle Diophantine inequalities. For example, inequalities $|F(\vec{n}) - a/q| < \delta$ where $F(\vec{n})$ is a linear form in irrational powers of $n_1, \dots, n_s \in \{1, \dots, N\}$, $a$ and $0 < q \leq Q$ are integers, and $\delta > 0$, govern the large values taken by exponential sums involving irrational powers of integers. Such inequalities are considered in \cite{Robert--Sargos_three-dimensional_exponential_sums_with_monomials}, and are related to problems involving Piatetski--Shapiro numbers $\lfloor n^c \rfloor$ \cite{Chapman--Chow--Holdridge_additive_Ramsey_theory_over_Piatetski-Shapiro_numbers}, and the Poissonian pair correlation property for the sequence $\alpha n^c$ \cite{Radziwill--Shubin}.
\end{enumerate}

\subsection{Notation and conventions}
\label{notations_and_conventions}

We use the Bachmann--Landau big-$O$ notation, so $f = O_{A, B, \dots}(g)$ for real-valued functions $f$ and $g$ if there exists a constant $c$, depending on the quantities $A, B, \dots$, such that $|f(x)| \leq c g(x)$ for a set of $x$ that will be clear from the context. We also use the Vinogradov notation, where $f \ll_{A, B, \dots} g$ if $f = O_{A, B, \dots}(g)$. We call the implicit constant $c$ in the big-$O$ or Vinogradov notation the big-$O$ constant. We stress that dependencies in the constant are always as written, so $f \ll g$ means the big-$O$ constant doesn't depend on any other parameters. We use the little-$o$ notation $f = o(g)$ to mean that $f(x)/g(x) \to 0$ as $x \to \infty$, and we use the notation $f \sim g$ to mean $f(x)/g(x) \to 1$ as $x \to \infty$. In the latter situation, one says that $f$ and $g$ are asymptotic. We also write $o_{A, B, \dots}$ and $\sim_{A, B, \dots}$ when the convergence is not necessarily uniform in the parameters $A, B$, etc. For example if $f(x) = x^2 + O_A(x)$ then $f \sim_A x^2$, and $O_A(x) = o_A(x^2)$.

%% file: article/content2/preliminaries.tex
\section{Preliminaries}
\label{preliminaries}

\subsection{Algebraic tori}

We recall facts and notions regarding algebraic tori. For details and proofs, one can consult \cite[Chapter 3]{Bombieri--Gubler_heights_in_Diophantine_geometry}. The algebraic torus is the algebraic group $\G_m$ over $\C$ whose group of $\C$-valued points is $\G_m(\C) = \C^\ast$. Thus the set $\R_{>0}$ of positive real numbers is a subgroup of $\G_m(\C)$. When we speak of the Zariski closure of a subset $X$ of $\R_{>0}^n$ in $\G_m^n$, we first identify $X$ with a subset of the closed points of $\G_m^n$, and then take the Zariski closure of the set of closed points.

A closed subvariety of $\G_m^n$ that is also an algebraic group with the same operation is called an algebraic subgroup. Thus $\R_{>0}$ is not an algebraic subgroup of $\G_m$, because it is not a closed subvariety. If $\Lambda$ is a subgroup of $\Z^n$, then the equations
\begin{align}
\label{monomial_relation}
    \vec{x}^{\vec{\lambda}} = x_1^{\lambda_1} \dots x_n^{\lambda_n} = 1
\end{align}
for each $\vec{\lambda} \in \Lambda$ define an algebraic subgroup. Conversely, every algebraic subgroup of $\G_m^n$ is given by multiplicative relations of this form, for some subgroup $\Lambda \subseteq \Z^n$. If $\Lambda$ has rank $r$, then the corresponding subgroup $H$ has dimension $n - r$, and $H$ is connected if and only if $\Lambda$ is a primitive lattice; that is, if and only if $(\Lambda \otimes_{\Z} \Q) \cap \Z^n = \Lambda$. If $H$ is a connected algebraic subgroup of dimension $k$, then $H$ is isomorphic to $\G_m^k$ as algebraic groups. 

Every morphism of algebraic groups $\phi : \G_m^n \to \G_m^k$ takes the form
\begin{align*}
    \phi(\vec{x}) = (\vec{x}^{\vec{a}_1}, \dots, \vec{x}^{\vec{a}_k})
\end{align*}
for some $\vec{a}_1, \dots, \vec{a}_k \in \Z^n$, where we have used, and will continue to use, the notation in (\ref{monomial_relation}). It is clear that the kernel of such a morphism of algebraic groups is the algebraic subgroup of $\G_m^n$ corresponding to the subgroup of $\Z^n$ generated by $\vec{a}_1, \dots, \vec{a}_n$. The morphisms of algebraic groups $\G_m^n \to \G_m$ are called the characters of $\G_m^n$. We shall also need the following fact about morphisms of algebraic tori.

\begin{lemma}
\label{preimage_of_orthant_under_monomial_map}
    Let $\phi : \G_m^n \to \G_m^k$ be an injective morphism of algebraic groups. If $\vec{x} \in \R^n$ and $\phi(\vec{x}) \in \R_{>0}^k$ then $\vec{x} \in \R_{>0}^n$.
\end{lemma}
\begin{proof}
    We can write $\vec{x} = \vec{y} \vec{z}$ under coordinate-wise multiplication, where $\vec{y} \in \{-1, 0, 1\}^n$ and $\vec{z} \in \R_{>0}^n$. Then $\phi(\vec{y})$ is in $\R_{>0}^k$ because $\phi(\vec{x})$ and $\phi(\vec{z})$ are. Now $\R_{>0}^k \cap \{-1, 0, 1\}^k$ is just the identity, so $\phi(\vec{y})$ is trivial, and so $\vec{y}$ is trivial, because $\phi$ is injective. The lemma follows.
\end{proof}

By abuse of notation, we will use $\exp$ to denote the usual exponential function $\exp : \C \to \C^\ast$, and its $n$-fold product $\exp : \C^n \to (\C^\ast)^n$.

\subsection{O-minimal geometry and point counting}

We first recall some aspects of semi-algebraic and o-minimal geometry. For more details, one can consult the book of Van den Dries \cite{Van_den_Dries_tame_topology_and_o-minimal_structures}. A semi-algebraic set $X$ in $\R^n$ is a finite union of sets of the form $\{\vec{x} \in \R^n : f(\vec{x}) = 0\}$ or $\{\vec{x} \in \R^n : f(\vec{x}) > 0\}$, where $f$ is a polynomial with real coefficients. A semi-algebraic set is called irreducible if its Zariski closure in affine space $\mathbb{A}^n$ is irreducible. We will be considering semi-algebraic subsets of $\R_{>0}^n$ and Zariski closures in $\G_m^n$. In this case, the Zariski closure of $X$ in $\G_m^n$ is irreducible if and only if the Zariski closure of $X$ in $\mathbb{A}^n$ is irreducible. The smooth locus of $X$ is exactly the points of $X$ corresponding to smooth closed points in the Zariski closure. If $X$ is irreducible, then its smooth locus is a real-analytic manifold. The dimension of $X$ is then equal to the dimension of its Zariski closure in $\G_m^n$, which is equal to its dimension as a real-analytic manifold. 

The fact that a semi-algebraic set remains semi-algebraic under a coordinate projection is usually attributed to Tarski \cite{Tarski_a_decision_method_for_elementary_algebra_and_geometry} and Seidenberg. This is a property known as quantifier elimination for the structure $\R_{\rm alg}$ of semi-algebraic sets, and it shows that $\R_{\rm alg}$ is a so-called o-minimal structure. Let us just say that a structure $S$ prescribes a collection of subsets of each Euclidean space $\R^n$, that are then said to be definable in the structure $S$, or $S$-definable. If the only $S$-definable subsets of $\R$ are finite unions of points and intervals, then the structure $S$ is said to be o-minimal.

We will also work in the structure $\R_{\exp}$ generated by the graph of the real exponential function. Thus a set $X \subseteq \R^n$ is $\R_{\exp}$-definable if it can be obtained from unions, intersections, complements, products, and coordinate projections of semi-algebraic sets or the graph of the exponential function. By a theorem of Wilkie \cite{Wilkie_model_completeness_results_for_expansions_of_the_ordered_field_of_real_numbers_by_restricted_Pfaffian_functions_and_the_exponential_functions}, the structure $\R_{\exp}$ is o-minimal.

The significance of o-minimal structures in our arguments lies in the Pila--Wilkie theorem and its variants. The multiplicative Weil height on the rational numbers is $H(a/q) = \max(|a|, q)$ for coprime integers $a$ and $q > 0$, and is extended to $\vec{a} \in \Q^n$ by $H(\vec{a}) = \max\{H(a_i) : i \in \{1, \dots, n\}\}$. For $X \subseteq \R^n$, definable in an o-minimal structure, we let $X^{\rm alg}$ denote the algebraic part of $X$, which is the union of all semi-algebraic subsets of $X$ having positive dimension, and we let $X^{\rm trans} = X \setminus X^{\rm alg}$ denote the transcendental part. We also let
\begin{align*}
    Y(\Q, N) = \{\vec{a} \in Y \cap \Q^n : H(\vec{a}) \leq N\}
\end{align*}
for $Y \subseteq \R^n$. Pila and Wilkie proved \cite[Theorem 1.8]{Pila--Wilkie_the_rational_points_on_a_definable_set} that if $X$ is definable in an o-minimal structure, then
\begin{align*}
    |X^{\rm trans}(\Q, N)| \ll_\eps N^\eps
\end{align*}
for every $N \geq 1$ and every $\eps > 0$. Let $X$ and $\Lambda$ be definable sets, and suppose there is a definable map $X \to \Lambda$, which means the graph of $X \to \Lambda$ is definable. In this case, one says that the fibres $X_\lambda$ of $X \to \Lambda$ form a definable family. In the same work, Pila and Wilkie prove a stronger result \cite[Theorem 1.9]{Pila--Wilkie_the_rational_points_on_a_definable_set}, that if $X \to \Lambda$ is a family that is definable in an o-minimal structure, then
\begin{align*}
    |X^{\rm trans}_\lambda(\Q, N)| \ll_{X, \eps} N^\eps
\end{align*}
for every $N \geq 1, \lambda \in \Lambda$, and every $\eps > 0$.

Wilkie conjectured (e.g., in \cite[Conjecture 1.11]{Pila--Wilkie_the_rational_points_on_a_definable_set}) that the $N^\eps$ appearing in the two results above could be improved to a power of $\log N$ when the o-minimal structure is $\R_{\exp}$. This conjecture has recently been settled by Binyamini, Novikov and Zak, and the following definable family version of their result is the main counting theorem we will use.

\begin{theorem}
[O-minimal point counting]
\label{bnz_counting_theorem}
    \cite[Equation (6)]{Binyamini--Novikov--Zak_Wilkie's_conjecture_for_Pfaffian_structures}
    Let $X \to \Lambda$ be an $\R_{\exp}$-definable family. There exists an effectively computable positive integer $\alpha$ such that
    \begin{align*}
        |X^{\rm trans}_\lambda(\Q, N)| \ll_X (\log N)^\alpha
    \end{align*}
    for every $N \geq 1$ and $\lambda \in \Lambda$.
\end{theorem}

%% file: article/content2/functional_transcendence.tex
\section{Functional transcendence}
\label{functional_transcendence_section}

In this section we prove our functional transcendence theorem (Theorem \ref{functional_transcendence_theorem}) for continuous group homomorphisms $\phi : \R_{>0}^n \to \R_{>0}^n$.

\begin{lemma}
    [Continuous group homomorphisms of $\R_{>0}^n$]
    \label{continuous_group_homomorphisms_of_the_orthant}
    Let $\phi : \R_{>0}^n \to \R_{>0}^n$ be a continuous group homomorphism. Then $\phi$ is real-analytic, and there are real vectors $\vec{a}_1, \dots, \vec{a}_n \in \R^n$ such that
    \begin{align*}
        \phi(\vec{x}) = 
        (\vec{x}^{\vec{a}_1}, \dots, \vec{x}^{\vec{a}_n})
    \end{align*}
    for all $\vec{x} \in \R_{>0}^n$, in the notation of (\ref{monomial_relation}).
\end{lemma}
\begin{proof}
    By pre-composing with the real exponential $\R^n \to \R_{>0}^n$ and post-composing by the real logarithm $\R_{>0}^n \to \R^n$ we obtain a continuous group homomorphism $\sigma : \R^n \to \R^n$. The group homomorphism property shows that $\sigma$ is $\Q$-linear and continuity shows that $\sigma$ is $\R$-linear. Thus $\sigma(\vec{x}) = A^t \vec{x}$ for some matrix $A$ with columns $\vec{a}_1, \dots, \vec{a}_n$. Therefore
    \begin{align*}
        \phi(\vec{x}) = \exp(A^t \log(\vec{x})) =
        (\vec{x}^{\vec{a}_1}, \dots, \vec{x}^{\vec{a}_n})
    \end{align*}
    as required.
\end{proof}

We now formulate a suitable transcendental or irrational property for endomorphisms of $\R_{>0}^n$.

\begin{definition}
    [Restricting to a morphism of algebraic groups]
    \label{restricting_to_a_morphism_of_algebraic_groups}
    Let $\phi : \R_{>0}^n \to \R_{>0}^n$ be a group homomorphism. We say that $\phi$ restricts to a morphism of algebraic groups, if there exists some algebraic subgroup $H$ of $\G_m^n$ having positive dimension, and a morphism $\sigma : H \to \G_m^n$ of algebraic groups, such that $\phi(\vec{x}) = \sigma(\vec{x})$ for all $\vec{x} \in \R^n_{>0} \cap H(\R)$; that is, if the restriction of $\phi$ to $H(\R)$ extends to a morphism on $H$.
\end{definition}

Note that $\phi$ restricts to a morphism of algebraic groups if and only if there is $\sigma : H \to \G_m^n$ as in Definition \ref{restricting_to_a_morphism_of_algebraic_groups} with $H$ of dimension $1$. Let us show that the condition in Definition \ref{restricting_to_a_morphism_of_algebraic_groups} is not a trivial one.

\begin{lemma}
    Let $H$ be an algebraic subgroup of $\G_m^n$. Then $\R_{>0}^n \cap H(\R)$ is a real-analytic Lie group with dimension equal to the dimension of $H$ as an algebraic subgroup of $\G_m^n$.
\end{lemma}
\begin{proof}
    Consider the exponential map $\exp : \C^n \to \G_m^n(\C)$. Since the exponential map is surjective, $\R_{>0}^n \cap H(\R) = \exp(\R^n \cap T)$, where $T$ is the complex tangent space of $H$. If $H$ corresponds to a subgroup $\Lambda \subseteq \Z^n$ as in Section \ref{preliminaries} then $T$ is the orthogonal complement of $\Lambda_\C$ in $\C^n$. In particular, $\R^n \cap T$ is the set of real points on a linear subvariety of $\C^n$ defined over the integers, and is therefore a real-analytic manifold of dimension equal to the dimension of $T$, which is in turn equal to the dimension of $H$ as an algebraic subgroup of $\G_m^n$. Since $\exp$ is a (real-analytic) diffeomorphism that is injective when restricted to $\R^n$, the image of $\R^n \cap T$ must also be a real-analytic manifold of the same dimension. Now a subgroup of a Lie group that is a manifold is again a Lie group.
\end{proof}

Consider $\phi(x, y) = (x^n, y^c)$, where $n$ is a non-zero integer and $c$ is real. Then $\phi$ restricts to the morphism $x \mapsto x^n$ on the algebraic subgroup defined by $y = 1$. Now consider $\phi(x, y) = (x^c, y^c)$, where $c$ is real and irrational. Then $\phi$ stabilises the algebraic subgroup defined by $x = y$, but it does not restrict to a morphism of algebraic groups there.

\begin{lemma}
[Criterion for restricting to a morphism of algebraic groups]
\label{restricting_to_a_morphism_of_algebraic_groups_criterion}
    Let $\phi$ be given as in Lemma \ref{continuous_group_homomorphisms_of_the_orthant}, and let $A$ be the matrix with columns $\vec{a}_1, \dots, \vec{a}_n$. Suppose $A^t \vec{b} = \vec{a}$ for some $\vec{a}, \vec{b} \in \Z^n$ with $\vec{b}$ non-zero. Then $\phi$ restricts to a morphism of algebraic groups.
\end{lemma}
\begin{proof}
    The action of $A^t$ on $\C^n$ restricts to the identity as a linear map between the $1$-dimensional linear subspaces $\C \vec{b}$ to $\C \vec{a}$. Since $\vec{b}$ and $\vec{a}$ are integer vectors, $\exp(\C \vec{b})$ and $\exp(\C \vec{a})$ are the complex points of two (not necessarily connected) algebraic subgroups $H_1$ and $H_2$. Moreover, if $\vec{x} \in \R_{>0}^n \cap H_1(\R)$, we can let $\vec{x} = \exp(u \vec{b})$ where $u \in \R$. Then $\phi(\vec{x}) = \exp(A^t(u\vec{b})) = \exp(u\vec{a}) = \vec{x}^{\vec{a}}$. Hence $\phi$ restricts to a morphism of algebraic groups. 
\end{proof}


\begin{theorem}
[Functional transcendence]
\label{functional_transcendence_theorem}
    Let $X$ be an irreducible, semi-algebraic subset of $\R_{>0}^n$. Let $\phi : \R_{>0}^n \to \R_{>0}^n$ be a continuous group homomorphism that does not restrict to a morphism of algebraic groups. Then the Zariski closure of $\phi(X)$ in $\G_m^n$ is a translate of a connected algebraic subgroup.
\end{theorem}
\begin{proof}
    The proof is by induction on $n$, the case $n = 0$ being trivial because then $X$ is a point. Thus assume $n > 0$ and that the theorem holds for algebraic tori of dimension smaller than $n$. If the theorem holds for $X$ then it holds for any translate of $X$ by a point of $\R_{>0}^n$. Thus we can assume that $X$ contains the identity $(1, \dots, 1)$ and that the identity is a smooth point of $X$. Let $d = \dim(X)$ be the dimension of the Zariski closure of $X$ in $\G_m^n$, which is equal to the dimension of $X$ as a real-analytic manifold. Since the exponential map is a diffeomorphism, the logarithm $\exp^{-1}(X)$ is again a real-analytic submanifold of $\R^n$ of dimension $d$. If $U \subseteq \exp^{-1}(X)$ is a connected, real-analytic chart containing $(0, \dots, 0)$, then there exists a parameterisation $(0, 1)^d \to U$ given by
    \begin{align*}
        \vec{y} \mapsto \vec{u}(\vec{y}) = (u_1(\vec{y}), \dots, u_n(\vec{y})),
    \end{align*}
    where $u_1, \dots, u_n$ are real-analytic functions vanishing at zero. In this case the map $\exp(\vec{u}(\vec{y})) : (0, 1)^d \to \R_{>0}^n$ parameterises a real-analytic chart of $X$ around $(1, \dots, 1)$.

    By Lemma \ref{continuous_group_homomorphisms_of_the_orthant}, $\phi$ is real-analytic; let $A$ be the matrix whose columns are the vectors $\vec{a}_1, \dots, \vec{a}_n$ in that lemma. Since the functions $\vec{u}$ vanish at zero, the functions $\vec{u}, A \vec{u}$ do not have constant terms in their power series expansions. Thus, by Ax's theorem \cite[Theorem 3]{Ax_on_Schanuel's_conjectures}, either
    \begin{align*}
        \trdeg_\C(\vec{u}, A\vec{u}, \exp(\vec{u}), \exp(A\vec{u})) \geq 2n + 
        \rank(J \, | \, A J),
    \end{align*}
    where $J_{ij} = \partial u_i/\partial y_j$ is the Jacobian matrix of the functions $\vec{u}$, or there is a non-trivial $\Z$-linear relation between $\vec{u}$ and $A\vec{u}$.

    Suppose the inequality holds. Since the functions $\vec{u}$ are components of a diffeomorphism, the rank of the matrix $(J \, | \, A J)$ is exactly $d$. Moreover, since $d$ is equal to the dimension of the Zariski closure of $X$ in $\G_m^n$, we have $\trdeg_\C(\exp(\vec{u})) = d$. Thus
    \begin{align*}
        2n + d \leq \trdeg_\C(\vec{u}) + 
        \trdeg_\C(\exp(\vec{u})) + 
        \trdeg_\C(\exp(A \vec{u}))
        \leq n + d + \trdeg_\C(\exp(A \vec{u}))
    \end{align*}
    by the trivial inequality $\trdeg_\C(\vec{u}) \leq n$. The map $\exp(A \vec{u})$ parameterises an open neighbourhood of the identity in the Zariski closure of $\phi(X)$ in $\G_m^n$, and the inequality above implies $\trdeg_\C(\exp(A\vec{u})) \geq n$. Thus the Zariski closure of $\phi(X)$ in $\G_m^n$ is $\G_m^n$ itself. This proves the theorem in this case.

    Thus assume that the alternative consequence of Ax's theorem holds, namely that there are $\vec{a}, \vec{b} \in \Z^n$, not both zero, such that
    \begin{align*}
        \vec{a} \cdot \vec{u} + \vec{b} \cdot (A\vec{u}) = (\vec{a} + A^t \vec{b}) \cdot \vec{u} = 0.
    \end{align*}
    If $\vec{a} + A^t \vec{b}$ vanishes then $A^t$ sends a non-zero integer vector to an integer vector, and so $\phi$ restricts to a morphism of algebraic groups, by Lemma \ref{restricting_to_a_morphism_of_algebraic_groups_criterion}. Thus $\vec{a} + A^t \vec{b}$ is non-zero, and we have found a non-trivial $\R$-linear dependence relation among the functions $\vec{u}$.
    
    Applying Ax's theorem \cite[Theorem 3]{Ax_on_Schanuel's_conjectures} to the functions $\vec{u}$ we have either
    \begin{align*}
        \trdeg_\C(\vec{u}, \exp(\vec{u})) \geq n + d
    \end{align*}
    or there is a non-trivial $\Z$-linear dependence relation between the functions $\vec{u}$. The inequality implies
    \begin{align*}
        n + d \leq \trdeg_\C(\vec{u}) + \trdeg_\C(\exp(\vec{u})) = 
        \trdeg_\C(\vec{u}) + d,
    \end{align*}
    which is inconsistent with the non-trivial $\R$-linear relation just constructed. Thus the $\R$-linear relation is promoted to a $\Z$-linear one: there exists some non-zero $\vec{c} \in \Z^n$ such that $\vec{c} \cdot \vec{u} = 0$. We can assume that $\vec{c}$ is has coprime coordinates. This shows that $\exp(U)$ is contained in a connected algebraic subgroup $H$ of codimension $1$. Since the Zariski closure of $\exp(U)$ is equal to the Zariski closure of $X$ (both in $\G_m^n$), we also see that $X$ is contained in $H$. It now suffices to apply the inductive hypothesis to $H$. Choose an isomorphism $\sigma : H \to \G_m^{n - 1}$ and let $Y = \sigma(X)$. Then $Y$ is irreducible, semi-algebraic, and $Y \subseteq \R_{>0}^{n - 1}$. Moreover, $\tau = \sigma \circ \phi \circ \sigma^{-1}$ is a continuous group homomorphism $\R_{>0}^{n - 1} \to \R_{>0}^{n - 1}$ by Lemma \ref{preimage_of_orthant_under_monomial_map}, that does not restrict to a morphism of algebraic groups because $\phi$ does not do so. By the inductive hypothesis $\overline{\tau(Y)}$ is a translate of an connected algebraic subgroup of $\G_m^{n - 1}$. Now since $\sigma$ is a homeomorphism for the Zariski topology
    \begin{align*}
        \sigma^{-1}(\overline{\tau(Y)}) =
        \overline{\phi(\sigma^{-1}(Y))} = \overline{\phi(X)}.
    \end{align*}
    Hence $\overline{\phi(X)}$ is a translate of a connected algebraic subgroup, and the theorem is proved.
\end{proof}

It is necessary to assume that $\phi$ does not restrict to a morphism of an algebraic subgroup of dimension at least $2$. Indeed, if $\phi$ restricts to a morphism $\G_m^2 \to \G_m^2$ given by $(x, y) \mapsto (x^a y^b, x^c y^d)$ where $a, b, c, d$ are integers, then the semi-algebraic curve
\begin{align*}
    x^a y^b + x^c y^d = 1
\end{align*}
in $\R_{>0}^2$ has Zariski closure $x + y = 1$, which is not a translate of an connected algebraic subgroup. Consider, however, $\phi(x, y) = (x^n, y^c)$ where $n$ is an integer and $c$ is a real irrational number. Then $\phi$ restricts to an endomorphism of the algebraic subgroup $y = 1$, so we cannot apply Theorem \ref{functional_transcendence_theorem}. Taking $X$ to be $x = y$, we find that $\phi(X)$ is Zariski dense in $\G_m^2$, so the consequence of the theorem still holds.

For our applications in Section \ref{proof_of_theorem_B}, we shall need the following.

\begin{lemma}
\label{lemma/raising_to_irrational_powers_not_a_morhpism_of_algebraic_subgroups}
    Let $c_1, \dots, c_n$ be real irrational numbers and let $\phi(\vec{x}) = (x_1^{c_1}, \dots, x_n^{c_n})$ for $\vec{x} \in \R_{>0}^n$. Then $\phi$ does not restrict to a morphism of algebraic subgroups.
\end{lemma}
\begin{proof}
    Suppose $\phi$ restricts to a morphism $\sigma : H \to \G_m^n$ defined on some non-trivial algebraic subgroup. Without loss of generality, $H$ corresponds to a lattice of rank $1$ in the tangent space of $\G_m^n$ generated by a vector $\vec{v}$. Then $\sigma$ is given by the action of integral matrix so it sends $\vec{v}$ to some integer vector. However this integer vector must agree with the image of $\vec{v}$ under the diagonal matrix with entries $c_1, \dots, c_n$. Since each $c_i$ are irrational, we must have $\vec{v} = 0$ so $H$ is trivial.
\end{proof}

%% file: article/content2/proof_of_theorem_B.tex
\section{Proof of Theorem \ref{additive_relations_in_irrational_powers}}
\label{proof_of_theorem_B}

In this section we shall apply the functional transcendence theorem (Theorem \ref{functional_transcendence_theorem}) to prove the following.

\begin{theorem}
[Additive relations in irrational powers]
    Let $r \leq s$ be positive integers. There exists an effectively constant $C_2 > 0$, depending on $s$, with the following property. Let $A$ be a subset of an $N$-term arithmetic progression in the non-negative reals. Let $c$ be a real irrational number. The number of solutions to the equation
    \begin{align}
    \label{additive_relations_in_irrational_powers_equation_proof}
        \sum_{i = 1}^s a_i^c = \sum_{j = 1}^r b_j^c
    \end{align}
    with $a_1, \dots, a_s, b_1, \dots, b_r \in A$ and $(a_1, \dots, a_s)$ not a permutation of $(b_1, \dots, b_r)$ is
    \begin{align*}
        O_s\left(
        |A|^{\theta}(\log N)^{C_2}
        \right),
    \end{align*}
    where $\theta = \max(1, \min(r, s - 1))$.
\end{theorem}

Since the equation (\ref{additive_relations_in_irrational_powers_equation_proof}) is invariant under dilations, we can assume $A$ lies in an arithmetic progression of the form $\{\alpha, \alpha + 1, \dots, \alpha + (N - 1)\}$. We first prove that one can take $\theta = \max(1, r)$ in Theorem \ref{additive_relations_in_irrational_powers}, which is the correct exponent when $s > r$.

\begin{lemma}
[Almost optimal bound when $s > r$]
\label{sum_of_irrational_powers_equals_fixed}
    For each fixed $v \in \R_{>0}$, the number of solutions to the equation
    \begin{align}
    \label{sum_of_irrational_powers_equals_fixed_equation}
        u_1^c + \dots + u_t^c = v
    \end{align}
    with $u_1, \dots, u_t \in A$ is $O_t((\log{N})^{O_t(1)})$. In particular the number of solutions to (\ref{additive_relations_in_irrational_powers_equation_proof}) with $a_1, \dots, a_s, b_1, \dots, b_r \in A$ and $(a_1, \dots, a_s)$ not a permutation of $(b_1, \dots, b_r)$ is
    \begin{align*}
        O_s(|A|^r(\log{N})^{O_s(1)}).
    \end{align*}
\end{lemma}
\begin{proof}
    Define the parameter space $\Lambda = \R \times \R_{>0} \times \R_{>0}$, where parameters $(c, \alpha, v)$ consist of the exponent $c$, the additive shift $\alpha$ defining the (dilated) arithmetic progression, and the constant $v$. Then the set
    \begin{align*}
        X = 
        \left\{
        (c, \alpha, v, \vec{x}) \in
        \Lambda \times \R_{>0}^t :
        v = (\alpha + x_1)^c + \dots + (\alpha + x_t)^c
        \right\},
    \end{align*}
    becomes an $\R_{\exp}$-definable family with the projection $X \to \Lambda$. Let $\lambda = (c, \alpha, v)$ with $c$ irrational. Let $Z$ be an irreducible semi-algebraic subset of $X_\lambda$ and let $W = t_{(\alpha, \dots, \alpha)}(Z)$, where $t_{(\alpha, \dots, \alpha)}$ is the additive translation by the vector $(\alpha, \dots, \alpha)$, which is also irreducible. Let $\phi : \R_{>0}^t \to \R_{>0}^t$ be given by $\phi(\vec{x}) = (x_1^c, \dots, x_t^c)$, so $\phi$ does not restrict to a morphism of algebraic groups by Lemma \ref{lemma/raising_to_irrational_powers_not_a_morhpism_of_algebraic_subgroups}. Then, by Theorem \ref{functional_transcendence_theorem}, the Zariski closure of $\phi(W)$ is a translate $\vec{g}H$ of an algebraic subgroup $H$ of $\G_m^t$, lying in the linear subvariety
    \begin{align*}
        v = x_1 + \dots + x_t.
    \end{align*}
    Let $\vec{g} = (g_1, \dots, g_t) \in \R_{>0}^t$. Then we have
    \begin{align*}
        v = g_1 \chi_1(h) + \dots + g_t \chi_t(h)
    \end{align*}
    for all $h \in H(\C)$, where $\chi_1, \dots, \chi_t$ are the coordinate characters of $\G_m^t$. By linear independence of characters and the fact that $\vec{g} \in \R_{>0}^t$, all $\chi_i$ must be trivial on $H$. Then $H$ is the trivial subgroup. Therefore $W$ and hence $Z$ are zero-dimensional. Thus $X_\lambda^{\rm alg} = \varnothing$ so the number of solutions to (\ref{sum_of_irrational_powers_equals_fixed_equation}) is
    \begin{align*}
        \ll_s |X_\lambda(\Q, N)| = |X_\lambda^{\rm trans}(\Q, N)| \ll_t
        (\log N)^{O_t(1)}
    \end{align*}
    by the o-minimal point counting result Theorem \ref{bnz_counting_theorem}. The second statement holds by applying the first statement to $v = b_1^c + \dots + b_r^c$ for the $|A|^r$ possible $(b_1, \dots, b_r)$.
\end{proof}

To prove Theorem \ref{additive_relations_in_irrational_powers}, it suffices to obtain the bound $O_s(|A|^{s - 1}(\log{N})^{O_s(1)})$ for $s > 1$, which is sharper only when $r = s$. Thus assume $r = s$. As in the proof of Lemma \ref{sum_of_irrational_powers_equals_fixed}, we will use a definable family. Define the parameter space $\Lambda = \R \times \R_{\geq 0} \times \R_{>0}^{s - 1}$, where parameters $(c, \alpha, \vec{a})$ correspond to the exponent $c$, the additive shift $\alpha$ defining the (dilated) arithmetic progression, and a vector $\vec{a} = (a_1, \dots, a_{s - 1})$ consisting of variables that we have fixed. Eventually we shall take $\vec{a} \in A^{s - 1}$. Now the set
\begin{align*}
    X = \left\{(c, \alpha, \vec{a}, x, \vec{y}) \in \Lambda \times \R_{>0}^{s + 1} :
    \sum_{i = 1}^{s - 1} a_i^c + (\alpha + x)^c = 
    \sum_{i = 1}^s (\alpha + y_i)^c\right\},
\end{align*}
becomes an $\R_{\exp}$-definable family when equipped with the projection $X \to \Lambda$. The following lemma classifies the semi-algebraic subsets of $X_\lambda$ when $c$ is irrational.

\begin{lemma}
\label{lemma/classification_of_semi-algebraic_curves}
    [Classification of semi-algebraic curves in $X_\lambda$]
    Let $\lambda = (c, \alpha, \vec{a}) \in \Lambda$ and suppose that $c$ is irrational. If $Z$ is an irreducible semi-algebraic curve in $X_\lambda$, then the Zariski closure of $Z$ has the form $t_{(\alpha, \dots, \alpha)}^*(\vec{g} H)$, where $H$ is a connected algebraic subgroup of $\G_m^{s + 1}$ of dimension $1$, $\vec{g} \in \R_{>0}^n$, and $t_{(\alpha, \dots, \alpha)}$ is the additive translation by the vector $(\alpha, \dots, \alpha)$. Moreover, there exists some $I \subseteq [s]$ with $1 \leq |I| \leq s - 1$ such that $H$ is defined by the equations $y_i = 1$ for $i \in I$ and $y_i = x$ for $i \in [s] \setminus I$, and such that the coordinates of the vector $\vec{g} = (g, g_1, \dots, g_s)$ satisfy the equations
    \begin{align*}
        \sum_{i \in I} g_i^c &= \sum_{i = 1}^{s - 1} a_i^c
        \\
        \sum_{i \in [s] \setminus I}
        g_i^c &= g^c.
    \end{align*}
\end{lemma}
\begin{proof}
    Let $Z$ be an irreducible semi-algebraic curve in $X_\lambda$, let $W = t_{(\alpha, \dots, \alpha)}(Z)$, and let
    \begin{align*}
        \phi(x, \vec{y}) = (x^c, y_1^c, \dots, y_s^c).
    \end{align*}
    Then, by Theorem \ref{functional_transcendence_theorem}, the Zariski closure of $\phi(W)$ is a translate $\vec{u} H$ of an connected algebraic subgroup $H$ of $\G_m^{s + 1}$ by some $\vec{u} \in \R_{>0}^{s + 1}$, lying in the linear subvariety
    \begin{align*}
        a_1^c + \dots + a_{s - 1}^c + x =
        y_1 + \dots + y_s.
    \end{align*}
    Let $\chi$ and $\chi_1, \dots, \chi_s$ denote the characters of $\G_m^{s + 1}$ corresponding to the projections to the coordinates $x, y_1, \dots, y_s$, respectively. If $\vec{u} = (u, u_1, \dots, u_s)$ then we have
    \begin{align*}
        a_1^c + \dots + a_{s - 1}^c + u\chi(h) =
        u_1\chi_1(h) + \dots + u_s\chi_s(h)
    \end{align*}
    for all $h \in H$. Suppose that $\chi$ is trivial on $H$. If each $\chi_i$ is also trivial on $H$ then $H$ is the trivial subgroup and $Z$ has dimension zero. Thus let $I \subseteq \{1, \dots, s\}$ be non-empty such that $i \in I$ if and only if $\chi_i$ is equal to a fixed, non-trivial character $\chi_I$ on $H$. By linear independence of characters we would then have
    \begin{align*}
        \sum_{i \in I} u_i = 0,
    \end{align*}
    which is absurd since $\vec{u} \in \R_{>0}^{s + 1}$. Therefore the character $\chi$ is non-trivial on $H$. By similar arguments, every character $\chi_1, \dots, \chi_s$ must be trivial or equal to $\chi$ on $H$. Let $I \subseteq [s]$ be such that $i \in I$ if and only if $\chi_i$ is trivial on $H$. Then $H$ is the $1$-dimensional subgroup defined by the equations $y_i = 1$ for all $i \in I$ and $y_i = x$ for all $i \in [s] \setminus I$, and linear independence of characters produces two relations
    \begin{align*}
        \sum_{i \in I} u_i &= \sum_{i = 1}^{s - 1} a_i^c, \\
        \sum_{i \in [s] \setminus I} 
        u_i &= u
    \end{align*}
    satisfied by the coordinates of $\vec{u}$. Now the Zariski closure of $\phi(W)$ is $\vec{u} H$ so $W$ is contained in $\phi^{-1}(\vec{u} H) = \phi^{-1}(\vec{u})H$. Let $\vec{g} = \phi^{-1}(\vec{u})$. Thus $Z$ is contained in $t_{(\alpha, \dots, \alpha)}^* (\vec{g}H)$ and if $\vec{g} = (g, g_1, \dots, g_s)$ then
    \begin{align*}
        \sum_{i \in I} g_i^c &= \sum_{i = 1}^{s - 1} a_i^c, \\
        \sum_{i \in [s] \setminus U} g_i^c &= g,
    \end{align*}
    as claimed. Since $Z$ and $t_{(\alpha, \dots, \alpha)}^* (\vec{g} H)$ are irreducible and of the same dimension, the Zariski closure of $Z$ must be as claimed. This completes the proof.
\end{proof}

We now finish the proof of Theorem \ref{additive_relations_in_irrational_powers}. Consider the parameter $\lambda = (c, \alpha, \vec{a})$ as fixed and suppose we have a solution to (\ref{additive_relations_in_irrational_powers_equation_proof}), which now consists of $(a_s, b_1, \dots, b_s)$, since $\vec{a} = (a_1, \dots, a_{s - 1})$ is fixed. Then $t_{(\alpha, \dots, \alpha)}^* (a_s, b_1, \dots, b_s)$ lies in $X_\lambda$. Suppose it lies on a semi-algebraic curve $Z$. Then $(a_s, b_1, \dots, b_s)$ lies in some translate $\vec{g}H$ of the form given in Lemma \ref{lemma/classification_of_semi-algebraic_curves}, parameterised by some $I \subseteq [s]$ with $1 \leq |I| \leq s - 1$. The relations satisfied by the coordinates of $\vec{g}$ then yield relations satisfied by $(a_1, \dots, a_s, b_1, \dots, b_s)$
\begin{align}
\label{equation/relations_for_semi-algebraic_solutions}
    \sum_{i = 1}^{s - 1} a_i^c &=
    \sum_{i \in I} b_i^c \\
    a_s^c &= \sum_{i \in [s] \setminus I} b_i^c.
\end{align}
If $r = s = 2$ then the relations take the form $a_1 - b_1 = a_2 - b_2 = 0$ or $a_1 - b_2 = a_2 - b_1 = 0$, so $(a_1, a_2)$ is a permutation of $(b_1, b_2)$. Thus, in this case, all non-trivial relations lie in $X_\lambda^{\textrm{trans}}$ and the total number of non-trivial relations is
\begin{align*}
    \sum_{a_1 \in A} |X_{(c, \alpha, a_1)}^{\textrm{trans}}(\Q, N)| \ll |A| (\log N)^{O(1)}.
\end{align*}
as claimed.

Now suppose $s$ is larger than $2$, and we have proved Theorem \ref{additive_relations_in_irrational_powers} for all smaller values of $s$. The number of solutions to the first equation in (\ref{equation/relations_for_semi-algebraic_solutions}) is $\ll_s |A|^{s - 2} (\log N)^{O_s(1)}$, unless $|I| = s - 1$ and $a_1, \dots, a_{s - 1}$ are a permutation of the $b_i$ for $i \in I$. In this case, the second equation in (\ref{equation/relations_for_semi-algebraic_solutions}) is simply $a_s^c = b_i^c$ for the unique $i \in [s] \setminus I$, so the coordinates of the vector $(a_1, \dots, a_s)$ are really a permutative of the coordinates of $(b_1, \dots, b_s)$. By an application of Lemma \ref{sum_of_irrational_powers_equals_fixed}, the second equation has $\ll_s |A|(\log N)^{O_s(1)}$ solutions. By the o-minimal counting theorem (Theorem \ref{bnz_counting_theorem}), the total number of solutions is
\begin{align*}
    \sum_{I \subseteq [r]}
    |A|^{s - 1}(\log{N})^{O_s(1)} +
    \sum_{\vec{a} \in A^{s - 1}}
    |X_{(c, \alpha, \vec{a})}^{\rm trans}(\Q, N)| \ll_s
    |A|^{s - 1}(\log{N})^{O_s(1)},
\end{align*}
where the first term counts the solutions that lie in $X_{(c, \alpha, \vec{a})}^{\rm alg}$ for some $\vec{a} \in A^{s - 1}$, and the second term counts the solutions lying in $X_{\lambda}^{\rm trans}$ for all $\lambda \in \Lambda$. This concludes the proof of Theorem \ref{additive_relations_in_irrational_powers}.

%% file: article/content2/proof_of_theorem_C.tex
\section{Proof of Theorem \ref{lower_bound_on_linear_form_in_irrational_powers}}
\label{proof_of_theorem_C}

In this section we prove the following theorem.

\begin{theorem}
[Non-vanishing for well-approximated exponents]
\label{lower_bound_on_linear_form_in_irrational_powers}
    Let $-B \leq a_1, \dots, a_s \leq B$ be integers, not all zero, and let $N$ be a positive integer. There exists an effectively computable, positive-valued function $\psi(B, N, s, q)$ such that the following holds. Suppose $c$ is an irrational real number, and there exist integers $a$ and $q > 0$ with
    \begin{align}
    \label{Diophantine_approximation_condition_proof}
        \frac{a}{q} < c \leq \frac{a}{q} + \psi(B, N, s, q).
    \end{align}
    Then $F(n_1, \dots, n_s) = a_1 n_1^c + \dots + a_s n_s^c$ is non-zero for all multiplicatively independent $n_1, \dots, n_s \in \{1, \dots, N\}$.
\end{theorem}
\begin{proof}
    Let $\eps = c - a/q$. We must show that if $\eps$ is chosen small enough in terms of $B, N, s, q$, then the form
    \begin{align*}
        F(\vec{n}) = a_1 n_1^c + \dots + a_s n_s^c
    \end{align*}
    does not vanish. Let $R(x, \eps)$ denote the remainder term such that
    \begin{align*}
        x^\eps = 1 + \eps\log{x} + R(x, \eps)
    \end{align*}
    and $R(x, \eps) \ll_N \eps^2$ for all $x \in [1, N]$. Writing $n_i^c = n_i^{a/q}n_i^\eps$ and expanding each term $n_i^\eps$ around zero gives
    \begin{align*}
        F(\vec{n}) &= \sum_{i = 1}^s a_in_i^{a/q} (1 + \eps \log{n_i} + R(n_i, \eps)) \\
        &=
        \alpha + \eps \Lambda + \sum_{i = 1}^s a_in_i^{a/q}R(n_i, \eps),
    \end{align*}
    where
    \begin{align*}
        \alpha = \sum_{i = 1}^s a_in_i^{a/q}, \quad
        \Lambda = \sum_{i = 1}^s a_in_i^{a/q}\log{n_i}.
    \end{align*}
    The idea of the proof is the following. If $\alpha$ does not vanish, then, being an algebraic integer with bounded conjugates, we can write down a lower bound for $\abs{\alpha}$. At the same time, we can make $\eps$ small enough that the lower bound on $\abs{\alpha}$ is much larger than the other terms contributing to $F(\vec{n})$. On the other hand, if $\alpha$ vanishes then a positive lower bound on $\Lambda$ can be used since the remainder term has order $\eps^2$. In both cases, we will need to use the estimate
    \begin{align*}
        \abs{\sum_{i = 1}^s a_i n_i^{a/q} R(n_i, \eps)} \ll_{B, N, s} \eps^2.
    \end{align*}
    First suppose $\alpha$ is non-zero. Using the reverse triangle inequality and the estimate $|\Lambda| \ll_{B, N, s} 1$ we have
    \begin{align*}
        \abs{F(\vec{n})} &\geq
        \abs{\alpha} - O_{B, N, s}(\eps).
    \end{align*}
    By the triangle inequality, the conjugates of $\alpha$ lie in the disc of radius $BN^cs$, and since $\alpha$ is an algebraic integer of degree at most $q^s$ we have
    \begin{align*}
        1 \leq \abs{\alpha}(BN^c s)^{q^s - 1}.
    \end{align*}
    In particular $|\alpha| \gg_{B, N, s, q} 1$. Hence, taking $\eps$ small in terms of $B, N, s, q$, we find that $F(\vec{n})$ is non-zero.

    Now suppose $\alpha$ is zero. In this case
    \begin{align*}
        |F(\vec{n})| \geq \eps(|\Lambda| - O_{B, N, s}(\eps)),
    \end{align*}
    and we can use an effective linear forms in logarithms result of Fel'dman \cite[Theorem 1]{Fel'dman_an_improvement_of_the_estimate_of_a_linear_form_in_the_logarithms_of_algebraic_numbers} that furnishes a lower bound $|\Lambda| \gg_{B, N, s, q} 1$. Again, $\eps$ may be chosen small enough in terms of $B, N, s, q$ so that $F(\vec{n})$ does not vanish.
    
    Finally, the quantity $\psi(B, N, s, q)$ is chosen small enough to be an admissible upper bound for $\eps$ in both cases. The required upper bounds on $\eps$ are effectively computable in both cases.
\end{proof}

\begin{lemma}
\label{existence_of_uncountably_many_c_and_effectively_computable_examples}
    Let $\psi_n(q)$ be a sequence of effectively computable, positive-valued, decreasing functions. There exist infinitely many effectively computable real numbers $c$ with the following property. For every $n \geq 1$ there exist coprime integers $a$ and $q > 0$ with
    \begin{align}
    \label{diophantine_approximation_property}
        \frac{a}{q} < c \leq \frac{a}{q} + \psi_n(q).
    \end{align}
\end{lemma}
\begin{proof}
    Let $(d_k)$ be a sequence of positive integers satisfying $d_k \geq k$ for every $k \geq 1$. Let $c$ be real and $a$ and $q > 0$ be coprime integers defined by
    \begin{align*}
        c = \sum_{k \geq 1} 2^{-d_k} = 
        \frac{a}{q} + \sum_{k > n} 2^{-d_k}.
    \end{align*}
    Then $q = 2^{d_n}$, $a/q < c$ and
    \begin{align*}
        c - \frac{a}{q} \leq \sum_{k \geq d_{n + 1}} 2^{-k} =
        2^{1 - d_{n + 1}}.
    \end{align*}
    If the sequence $(d_k)$ has 
    \begin{align}
    \label{condition_required_for_binary_expansion}
        2^{d_{k + 1}} \geq \frac{2}{\psi_k(2^{d_k})} + e_k
    \end{align}
    for some sequence $(e_k)$ of non-negative integers and for all $k \geq 1$, then $2^{1 - d_{n + 1}} \leq \psi_n(q)$, as required. Finally, it is clear that if $\psi_n(q)$ is effectively computable, then an effectively computable sequence $(e_k)$ gives rise to an effectively computable sequence $(d_k)$ defined by letting $d_{k + 1}$ be minimal such that (\ref{condition_required_for_binary_expansion}) holds.
\end{proof}

To apply Lemma \ref{existence_of_uncountably_many_c_and_effectively_computable_examples} to the function $\psi(B, N, s, q)$ in Theorem \ref{lower_bound_on_linear_form_in_irrational_powers}, one considers the sequence of functions
\begin{align*}
    \psi_n(q) = \max \{\psi(B, N, s, q) : 1 \leq B, N, s \leq n\}.
\end{align*}
This proves the assertion in the introduction that if $S$ is a set of multiplicatively independent integers, then $S^{[c]}$ is linearly independent over $\Q$ for infinitely many effectively computable numbers $c$.

%% file: article/content2/rational_exponents.tex
\section{Rational exponents}
\label{rational_exponents}

In this section we consider the $2$-fold additive energy of $A^{[c]}$ when $A = \{1, \dots, N\}$ and $c$ is rational, with the objective of finishing the deduction of Corollary \ref{expansion_for_all_exponents}, which states that $|A^{[c]} + A^{[c]}| \sim N^2/2$ as $N \to \infty$ for all $c \in \R \setminus \{0, 1, 2\}$. Theorem \ref{expansion} covers the case when $c$ is irrational.

When $c$ is a positive integer, the question of solutions to the equation
\begin{align}
\label{energy_equation_rational_exponents}
    n_1^c + n_2^c = n_3^c + n_4^c
\end{align}
with $n_1, n_2, n_3, n_4 \in \{1, \dots, N\}$ has received plenty of attention over the years. We have already mentioned Landau's result \cite{Landau_on_the_partition_of_positive_integers_in_four_classes_according_to_the_minimal_number_of_squares_needed_to_their_additive_composition} that $E_2(A^{[2]}) \gg N^2(\log{N})^{1/2}$. Hooley addressed the case $c = 3$ in \cite{Hooley_on_the_representations_of_a_number_as_a_sum_of_two_cubes, Hooley_on_the_numbers_that_are_representable_as_the_sum_of_two_cubes} and general positive integral $c$ in many articles, e.g., \cite{Hooley_on_the_representation_of_a_number_as_the_sum_of_two_h-th_powers, Hooley_on_another_sieve_method_and_the_numbers_that_are_a_sum_of_two_h-th_powers}---in the second article, the conjecture that $A^{[c]}$ is a Sidon set when $c \geq 5$ is a positive integer is stated. For $c = 4$, the problem was considered by Greaves \cite{Greaves_on_the_representation_of_a_number_as_a_sum_of_two_fourth_powers} and an improvement was obtained for $c \geq 5$ by Skinner--Wooley \cite{Skinner--Wooley_sums_of_two_kth_powers}, utilising uniform bounds for integral points on curves, proved by Bombieri and Pila using the determinant method \cite{Bombieri--Pila}. The development of the $p$-adic determinant method by Heath-Brown and others yielded further improvements \cite{Heath-Brown_the_density_of_rational_points_on_curves_and_surfaces, Browning_equal_sums_of_two_kth_powers}. Recently, these methods have been pushed further by Salberger \cite{Salberger_counting_rational_points_on_projective_varieties}.

We note that the expectation for the number of non-trivial solutions can be explained in terms of the geometry of the surface $X_c$ in $\P^3$ defined by the equation (\ref{energy_equation_rational_exponents}). The surface $X_3$ is rational, and so its rational points are Zariski dense. The surface $X_4$ is K3 and admits an elliptic fibration over $\Q$ with generic fibre having Mordell--Weil rank $1$. The rational points are Zariski dense, and even dense in the real topology, by a result of Swinnerton-Dyer \cite{Swinnerton-Dyer_A4+B4=C4+D4_revisted}. Finally, the surfaces $X_c$ for $c \geq 5$ are of general type, and the conjecture that (\ref{energy_equation_rational_exponents}) admits no non-trivial solution when $c \geq 5$ can be regarded as a refinement of the Bombieri--Lang conjecture, wherein the rational points of $X_c$ are not Zariski dense.

The case of negative integers $c$ can be reduced to the case of positive $c$ by multiplying by $(n_1n_2n_3n_4)^{-c}$ to obtain a singular surface of degree $-3c$. The case $c = -1$ was studied in detail by Heath-Brown \cite{Heath-Brown_Cayleys_cubic}, where he obtained the order of magnitude $N(\log{N})^7$ for the number of non-trivial solutions, as predicted by Manin's conjecture.

For our purposes, we only need the following.

\begin{lemma}
[Power saving error term for integer exponents]
\label{power_saving_for_integer_exponents}
    Let $c \in \Z \setminus \{0, 1, 2\}$ and let $A = \{1, \dots, N\}$. Then there exists $\delta$, depending on $c$, such that
    \begin{align*}
        E_2(A^{[c]}) = 2N^2 + O_c(N^{2 - \delta}).
    \end{align*}
\end{lemma}
\begin{proof}
    When $c \geq 2$, we can use \cite[Corollary 0.7]{Salberger_counting_rational_points_on_projective_varieties}. When $c < 0$, we multiply by $(n_1n_2n_3n_4)^c$ to obtain
    \begin{align*}
        (n_2n_3n_4)^c + (n_1n_3n_4)^c =
        (n_1n_2n_4)^c + (n_1n_2n_3)^c
    \end{align*}
    and then apply \cite[Theorem 0.5]{Salberger_counting_rational_points_on_projective_varieties}.
\end{proof}

It will be convenient to introduce the notation $B(c, N)$ for the number of non-trivial solutions to (\ref{energy_equation_rational_exponents}) in $\{1, \dots, N\}$. Thus it remains to show $B(c, N) = o(N^2)$ for all rational $c$ that are not integers.

\begin{lemma}
[Reduction to integer exponents]
\label{reduction_to_integer_case}
    Let $a$ and $q > 0$ be coprime integers such that $c = a/q$. Then
    \begin{align*}
        B(a/q, N) = \sum_{b \leq N}
        B(a, (N/b)^{1/q}).
    \end{align*}
\end{lemma}
\begin{proof}
    Suppose $(n_1, n_2, n_3, n_4)$ is a non-trivial solution of (\ref{energy_equation_rational_exponents}). Let $n_i = a_i^q b_i$ where $a_i^q, b_i \in \{1, \dots, N\}$, and each $b_i$ is $q$-th power free. Consider the set $S = \{b_1^{a/q}, b_2^{a/q}, b_3^{a/q}, b_4^{a/q}\}$. An elementary argument shows that elements of $S$ are pairwise coprime. By \cite[Theorem 1.1]{Carr--O'Sullivan_on_the_linear_independence_of_roots}, the set $S$ is linearly independent over $\Q$, despite the fact that
    \begin{align*}
        a_1^a b_1^{a/q} + a_2^a b_2^{a/q} = a_3^a b_3^{a/q} + a_4^a b_4^{a/q}.
    \end{align*}
    Since the solution is non-trivial, we must have all $b_i$ equal, say $b_i = b$. Then $(a_1, a_2, a_3, a_4)$ is a solution of (\ref{energy_equation_rational_exponents}) with $c = a$ and each $a_i^q \leq N/b$. Hence the lemma.
\end{proof}

\noindent\textbf{Proof of Corollary \ref{expansion_for_all_exponents}.}
We use the power saving $B(a, M) \ll M^{2 - \delta}$ for integer exponents $a \not \in \{0, 1, 2\}$ in Lemma \ref{power_saving_for_integer_exponents}, together with the elementary estimate
\begin{align*}
    \sum_{b \leq N} \left(\frac{N}{b}\right)^{2 - \delta} \ll_\delta
    N(\log{N}) + N^{2 - \delta}
\end{align*}
for all $\delta > 0$. By Lemma \ref{reduction_to_integer_case}, this yields the corollary unless $c$ takes the form $1/q$ or $2/q$. When $a = 1$ we must have $q \geq 2$ so
\begin{align*}
    B(1/q, N) \ll_q \sum_{b \leq N} \left(\frac{N}{b}\right)^{3/2}
    \ll_q N^{3/2},
\end{align*}
and when $a = 2$ we must have $q \geq 3$ so
\begin{align*}
    B(2/q, N) \ll_q \sum_{b \leq N} \frac{N}{b} \ll_q N(\log{N}).
\end{align*}
This shows that $B(c, N) = o(N^2)$ for all $c \in \Q \setminus \{0, 1, 2\}$, and finishes the proof of Corollary \ref{expansion_for_all_exponents}.

%% file: article/backmatter/bibliography.tex

\bibliographystyle{alpha}
\bibliography{article/backmatter/references}